%% file: virtual_betti10_26_05.tex
\documentclass{amsart}
\usepackage{amsmath,amsthm,amscd,amssymb}
\usepackage{pb-diagram,enumerate}
\usepackage{float,epsfig}
\usepackage{graphicx}

\input amsmath-defs.tex

\newtheorem{thm}{Theorem}[section]
\newtheorem{prop}[thm]{Proposition}
\newtheorem{lem}[thm]{Lemma}
\newtheorem{example}{Example}
\newtheorem{cor}[thm]{Corollary}
\theoremstyle{definition}
\newtheorem*{defn}{Definition}

\newtheorem*{VHC}{Waldhausen's Virtual Haken Conjecture (VHC)}
\newtheorem*{VPBNC}{Virtual Positive Betti Number Conjecture (VPBNC)}
\newtheorem*{VIBNC}{Virtual Infinite Betti Number Conjecture (VIBNC)}
\newtheorem*{VFC}{Thurston's Virtual Fibering Conjecture (VFC)}
\newtheorem*{A}{Question A}
\newtheorem*{1}{Theorem}

\newcommand{\image}{\operatorname{image}}
\newcommand{\ke}{\operatorname{ker}}
\newcommand{\cok}{\operatorname{cokernel}}
\newcommand{\intt}{\operatorname{int}}
\newcommand{\tor}{\operatorname{Tor}}
\newcommand{\lk}{\ensuremath{\ell k}}

\begin{document}

\title{The Growth Rate of the First Betti Number in Abelian Covers of $3$-Manifolds }

\author{Tim D. Cochran}
\address{Rice University, Houston, Texas, 77005-1892}
\email{cochran@math.rice.edu}
\author{Joseph Masters}
\address{SUNY at Buffalo, Buffalo,N.Y.}
\email{jdmaster@eng.buffalo.edu}
\thanks{The first author was partially supported by the National Science Foundation. The second author was
supported by a National Science Foundation Postdoctoral Research Fellowship}

\maketitle

\section*{Abstract} 
 We give examples of closed hyperbolic 3-manifolds
 with first Betti number $2$ and $3$
for which no sequence of finite abelian covering spaces increases the
 first Betti number. For $3$-manifolds $M$ with first Betti number $2$
 we give a characterization in terms of some generalized self-linking numbers
 of $M$, for there to exist a family of $\mathbb{Z}_n$ covering spaces,
 $M_n$, in which $\beta _1(M_n)$ increases linearly with $n$.
 The latter generalizes work of M. Katz and C. Lescop [KL],
 by showing that the non-vanishing of any one of these invariants
 of $M$ is sufficient to guarantee certain optimal systolic inequalities for $M$ (by work of Ivanov and Katz
[IK]).

\section*{Introduction} Motivated by Waldhausen's work on Haken manifolds,
 and by W. Thurston's {\bf Geometrization Conjecture},
 it has been variously conjectured that, if $M$ is an orientable,
irreducible closed $3$-manifold with infinite fundamental group, then:

\begin{VHC} $M$ is finitely covered by a Haken manifold;
\end{VHC}

\begin{VPBNC} Some finite cover of $M$ has positive first Betti
number;
\end{VPBNC}

\begin{VIBNC} Either $\pi_1(M)$ is virtually solvable or $M$ has
finite covers with arbitrarily large first Betti number;
\end{VIBNC}
\begin{VFC} $M$ has a finite cover that fibers over the circle.
\end{VFC}
There are easy implications VIBNC$\Longrightarrow$VPBNC$\Longrightarrow$VHC
and VFC$\Longrightarrow$VPBNC $\Longrightarrow$VHC. 
Each implies, if $M$ is atoroidal, the
long-standing conjecture of Thurston that such a manifold admits a
geometric structure. It is interesting to note that even if $M$ is {\bf
assumed} to be hyperbolic, the conjectures above are open.

In this paper, we restrict our attention to VIBNC. (We note in passing
that the alternative ``$\pi_1(M)$ is virtually solvable'' is sometimes
replaced by the a priori stronger alternative that ``$M$ is finitely
covered by the 3-torus, a nilmanifold or a solvmanifold.'') One rich
source of finite covering spaces is the set of iterated (regular)
finite {\bf abelian} covering spaces. Thus specifically, in this paper we
consider the question:

\begin{A} Does there exist
 an integer $m$, such that, if $M$ is any closed, atoroidal $3$-manifold
 with $\beta_1(M) \geq m$ then  $\b_1(M)$ can be increased
in a finite abelian covering space?
\end{A}

 Note that some condition
 on $H_1(M)$ is necessary, for if $H_1(M)=0$, then $M$ admits no
non-trivial abelian covering spaces. Counter-examples also exist for
many manifolds with $\b_1(M)=1$. For if $M$ is zero-framed surgery on a
knot in $S^3$, then it is easy to show that $H_1(\wt M;\BQ)\cong\BQ\op
Q[t,t^{-1}]/\<\Delta_k,t^n-1\>$ where $\wt M$ is the $n$-fold cyclic
cover and $\Delta_k$ is the Alexander polynomial of $K$. Thus
$\b_1(\wt M)=\b_1(M)=1$ except when $\Delta_k$ has a cyclotomic factor.
We begin this paper by observing that counter-examples also exist in the
cases $\b_1(M)=2$ and $\b_1(M)=3$. 

\begin{1}\label{mainthm1} \textit{There exist closed hyperbolic $3$-manifolds $M$
with $\b_1(M)=2$ (respectively $3$) for which no sequence of finite abelian
covers increases the first Betti number.
 More generally, if a sequence of regular covers of M increases
 the first Betti number, then one of the covering groups contains a
 non-trivial perfect subgroup.}
\end{1}

It is noteworthy that Question A is still open.

 If $\b_1 > 0$, then there is an epimorphism
 $\pi_1(M)\to \mathbb{Z}$, and a corresponding sequence of finite cyclic
 covers of $M$.
 Our second contribution is, in the case $\b_1(M)=2$,
 to give necessary and sufficient conditions,
 of a somewhat geometric flavor, for the Betti number
 of these covers to increase linearly with the covering degree.
 This is the content of Section 2.

\section{On abelian covers of hyperbolic 3-manifolds
 with $\b_1(M)=2$ and $3$}\label{failure} In this section,
 we observe that, if Question A has an affirmative answer,
 then the integer $m$ must be at least 4.

\begin{thm}\label{failure23} There exist closed hyperbolic $3$-manifolds
$M$ with $\b_1(M)=2$ (respectively $3$) such that if $\wt M$ is obtained
from $M$ by taking a sequence of finite abelian covering spaces, then
$\b_1(\wt M)=\b_1(M)$.
 More generally, if a sequence of regular covers of M increases
 the first Betti number, then one of the covering groups contains a
 non-trivial perfect subgroup.
\end{thm}

\begin{proof} Begin with a ``seed'' manifold $N$ whose fundamental group
is nilpotent.  Recall that the {\bf Heisenberg manifold} with Euler class $e$ is the circle bundle over the
torus with Euler class $e$. The fundamental group of such a $3$-manifold is the nilpotent group $\langle x,y,t
: [x,y]=t^e, [x,t], [y,t]\rangle$, called the {\bf Heisenberg group} of Euler class $e$. For our seed manifold
with $\b_1(N)=2$, we shall take $N$ to be the Heisenberg manifold with Euler class $1$, that can also be
described as $0$-framed surgery on a
Whitehead link. Thus in this case $\pi_1(N)\cong F/F_3$ where $F$ is the free group of rank~2 and $F_3$ is the
third term of the
lower-central series of $F$. When $\b_1(N)=3$, we take our seed manifold $N$ to be $S^1\x S^1\x S^1$, the
Heisenberg manifold of Euler class $0$. Note that each of the Heisenberg groups of non-zero Euler class has
$\beta_1=2$ while the Heisenberg group of Euler class $0$ has $\beta_1=3$.
 
First we claim that \emph{no} finite cover of $N$ will increase the first Betti number, which follows
immediately from the Lemma below (which surely is well-known to experts).
\begin{lem}\label{nilpotent} Suppose $A$ is a Heisenberg group with non-zero (respectively zero) Euler class.
If $\wt A$ is any finite index subgroup of $A$, then $\wt A$  is a Heisenberg group of non-zero (respectively,
zero) Euler class. Hence in all cases $\beta_1(\wt A) = \beta_1(A)$.
\end{lem}
\begin{proof}[Proof of Lemma~\ref{nilpotent}] The result is obvious for $A=\mathbb{Z}\x \mathbb{Z}\x
\mathbb{Z}$ so we assume that $A$ is a Heisenberg group of non-zero Euler class $e$. Then $A$ is a central
extension as shown below.
$$
1\lra \mathbb{Z}\overset{i}{\lra} A \overset{\pi}{\lra} \mathbb{Z}\times\mathbb{Z}\lra 1
$$
Since $\wt A$ is a finite index subgroup of $A$, $\pi(\wt A)$ is a finite index subgroup of
$\mathbb{Z}\times\mathbb{Z}$ which is hence isomorphic to $\mathbb{Z}\times\mathbb{Z}$. Moreover the kernel of
the map $\pi:\wt A\to \pi(\wt A)$ is a finite index subgroup of kernel($\pi$)$=\mathbb{Z}$ which is contained
in the center of $A$. It follows that $\wt A$ is also a central extension of the above form  and hence is also
a Heisenberg group. We claim that $\wt A$ has non-zero Euler class. Suppose not. Then $\wt A$ is abelian. But
$A\cong \langle x,y,t : [x,y]=t^e, [x,t], [y,t]\rangle$ where $e\neq 0$. Consider the elements $\{x,y\}$. There
is some positive integer $n$ such that both $x^n$ and $y^n$ lie in the subgroup $\wt A$ where they commute.
Thus $[x^n,y^n]=1$ in $A$. However since $[x,y]=t^e$, and $t$ commutes with $x$ and $y$, it is easy to see that
$x^ny^n=t^ky^nx^n$ where $k=n^2e$ and so $1=[x^n,y^n]=t^k$. This implies that $t$ is of finite order.  
 However, any Heisenberg group is the fundamental group of a circle bundle over the torus, which is
an aspherical 3-manifold. Thus $A$ has geometric dimension $3$ and cannot have torsion, for a contradiction.
\end{proof}

Next alter the seed manifold in a subtle way using the following result of
A. Kawauchi [Ka1 p. 450-452 , and Ka2 Corollary 4.3] (see also Boileau-Wang
 [BW section 4]).

\begin{prop}\label{hyperbolic} (Kawauchi) For any closed $3$-manifold $N$,
there exists a hyperbolic $3$-manifold $M$ and a degree 1 map $f:M\to N$
that induces an isomorphism on homology groups with local coefficients in
$\pi_1(N)$. Equivalently, if $\wt N$
 is  any covering space of $N$
 and $\wt f:\wt M  \to \wt N$ is the pull-back, then
 $\wt f$ induces isomorphisms on homology groups.
\end{prop}

\begin{proof} To the best of our knowledge, this result was first established by Kawauchi using his theory of
{\bf almost
identical imitations}. We sketch a proof using the approach of Boileau and Wang (which overlaps substantially with Kawauchi's approach). Recall that any $3$-manifold
$N$ contains a knot $J$ whose exterior is
hyperbolic. With more work, Boileau and Wang ensure that there exists
 such a knot
$J$ which is ``totally null-homotopic'', i.e., bounds a map of a 2-disk,
$\phi:D^2\to N$, such that the inclusion map $\pi_1(\image\phi)\to\pi_1(N)$
is trivial. Let $M_n$ be the result of $1/n$-Dehn surgery
on $N$ along $J$. By work of W. Thurston, for almost all
  $n$, $M_n$ is hyperbolic. Choose such an $M_n$ and denote it by $M$.
 Since $J$ is null-homotopic there is
a degree one map $f:M\to N$ that induces an isomorphism on $H_1$.

 Let $\wt N$ be a cover of $N$.  Since
 $J$ is null-homotopic, it lifts to $\wt N$, and there
 is an induced cover $\wt M$ and an induced map
 $\tl f:\wt M\to\wt N$.
 Since $J$ is totally null-homotopic, the pre-images
 of $J$ bound disjoint Seifert surfaces in $\wt M$,
 and so  $\tl f:\wt M\to\wt N$ is an
 isomorphism on homology.

\end{proof}

 For any map $f: M \to N$ satisfying the conclusion
 of Proposition \ref{hyperbolic}, $ker(f_*)$ is a perfect group.
 Indeed, Proposition~\ref{hyperbolic} states that for {\bf any}
covering space $\wt M$ of $M$ that is ``induced''
 from a cover $\wt N$ of $N$, the induced map $\tl f:\wt M\to\wt N$ is an
isomorphism on homology, so $\b_1(\wt M)=\b_1(\wt N)$. Specifically, letting $\wt N$ be the universal cover, 
$H_1(\wt M)\cong H_1(\wt N)=0$ showing that that $\pi_1(\wt M)$ is a perfect group. But $\pi_1(\wt M)$ is
kernel$(f_*)$. (Indeed, the condition that $f:M\to N$ induce an isomorphism on first homology with local
coefficients in
$\pi_1(N)$ is equivalent to the condition that the kernel of $f_*:\pi_1(M)\to\pi_1(N)$ be a perfect group).

Returning to the proof of our theorem, recall that $N$ is our seed
 Heisenberg manifold, and let $M$ be the manifold guaranteed by Proposition
 \ref{hyperbolic}.
 We claim that
 the manifold $M$ satisfies the conclusion of the theorem. For suppose
$\wt M\overset{p}{\lra}M$ is a regular finite covering space
 of $M$ corresponding to a surjection $\psi:\pi_1(M)\to F$, where $F$ is a
finite group that contains no nontrivial perfect subgroup (for example if $F$ is abelian). Then, since the
kernel of $f_*:\pi_1(M)\to\pi_1(N)$ is a perfect group $P$, and the perfect subgroup $\psi(P)\subset F$ must be trivial,
$\psi$ factors through $f_*:\pi_1(M)\to\pi_1(N)$ via a surjection
$\phi:\pi_1(N)\to F$. Therefore there is a finite regular cover $\wt N$ of
$N$ and a lift $\tl f:\wt M\to\wt N$. Notice that the only property of $M$ and $N$ needed for this argument is
that the kernel of $f_*:\pi_1(M)\to\pi_1(N)$ is a perfect group. 
Proceeding, by Proposition~\ref{hyperbolic} $H_1(\wt M)\cong H_1(\wt N)$ and by Lemma~\ref{nilpotent},
$\beta_1(\wt N)=\beta_1(N)$. Since $\beta_1(M)=\beta_1(N)$ we conclude that $\beta_1(\wt M)=\beta_1(M)$. This
shows that the first Betti number of $M$ cannot
 be increased by a \emph{single} regular $F$-cover unless $F$ contains a nontrivial perfect subgroup. In
particular, it shows that the first Betti number of $M$ cannot
 be increased by a single \emph{abelian} cover.

 Now suppose that $M_k \to ... \to M_0 = M$ is a sequence
 of regular covers, with covering groups $F_1, ..., F_k$,
 where no $F_i$ contains a nontrivial perfect subgroup.
 In the last paragraph we showed that the cover $M_1 \to M_0$ is the pull-back of a corresponding cover $N_1
\to N_0$.
 We claim that the kernel, $P_1$, of the lift $(f_1)*:\pi_1(M_1)\to\pi_1(N_1)$ is \emph{equal to} the kernel,
$P_0$, of $f_*:\pi_1(M_0)\to\pi_1(N_0)$ (here we view $\pi_1(M_1)$ as a subgroup of $\pi_1(M_0)$). For,
obviously $P_1 \subset P_0$
 and since $F_0$ contains no perfect subgroups, $P_0\subset P_1$. Thus $P_1$
 is a perfect group and thus $\wt f_1$ induces an isomorphism
 on homology (even with twisted coefficients).
 Thus we have recovered the inductive hypothesis of the previous
 paragraph and continuing
inductively, we get a sequence of finite covers
 $N_k \to ... \to N_1$, with $\b_1(M_k) = \b_1(N_k)$.
 Therefore, to finish the proof we only need to observe that $\b_1(N)$ cannot
 be increased by any sequence of finite covers, which was shown in Lemma~\ref{nilpotent}. 
 
\end{proof}

\section{Linear Growth of Betti Numbers in Cyclic Covering
Spaces}\label{lineargrowth}

In this section we ask whether or not it is possible to increase the first Betti number with \emph{linear growth rate} in some \emph{compatible family} of cyclic covering spaces. If $M_\infty$ is a fixed infinite cyclic covering space corresponding to an epimorphism $\psi :\pi_1(M)\to \mathbb{Z}$ then by a \emph{compatible family} we mean the usual family of finite cyclic covers $M_n$ associated to $\pi_1(M)\to \mathbb{Z}\to \mathbb{Z}_n$. By a \emph{linear growth rate} we mean $\varinjlim (\beta_1(M_n)/n)$ is positive. It was already known that linear growth occurs precisely when $H_1(M_\infty)$ has positive rank as a $\mathbb{Z}[t,t^{-1}]$-module \cite[Theorem 0.1]{Lu2}\cite[pg.35 Lemma 1.34,pg.453]{Lu1}. Therefore our contribution is to offer a more geometric way of viewing this criterion. We also point out an application to certain optimal systolic inequalities for such $3$-manifolds as have appeared in work of Katz [IK][KL].

One should note from the outset that if $\pi_1(M)$ admits an epimorphism to $\mathbb{Z}\ast\mathbb{Z}$, then it is an easy exercise to show that $\beta_1(M)$ can be increased linearly in finite cyclic covers since the same is patently true of the wedge of two circles. Such manifolds arise, for example, as $0$-framed surgery on $2$-component boundary links. This condition is not necessary, however, as we shall see in Example~\ref{example3} below.

Suppose $M$ is a closed, oriented $3$-manifold with $\b_1(M)=2$. Given any
basis $\{x,y\}$ of $H^1(M,\BZ)$ we shall define a sequence of higher-order
invariants $\b^n(x,y)$; $n\ge1$ taking values in sets of rational numbers.
 The invariants can be interpreted as certain Massey products in $M$. The
invariant $\b^1(x,y)$ is always defined, is independent of basis, and
essentially coincides with the invariant $\la$, an extension of Casson's
invariant, due to Christine Lescop [Les].
 If $\b^i$ is defined for all
$i<n$ and is zero, then $\b^n$ is defined (this is why the invariants are
called higher-order). If $H_1(M)$ has no torsion, so that $M$ can be
viewed as $0$-framed surgery on a 2-component link in a homology sphere
(with Seifert surfaces dual to $\{x,y\}$) then $\b^n$, when defined, is
the same as the sequence of link concordance invariants of the same name
due to the first author [C1]. In this case $\b^1$ was previously known as
the Sato-Levine invariant. 

After defining the invariants $\b^n(x,y)$, we show that their vanishing is equivalent to the linear growth of Betti numbers in the family corresponding to the infinite cyclic cover associated to $x$.

\begin{thm}\label{linear} Let $M$ be a closed oriented $3$-manifold with
$\b_1(M)=2$. The following are equivalent.
\begin{enumerate}
\item[A.] There exists a compatible family $\{M_n|n\ge1\}$ of finite
cyclic covers of $M$ such that $\b_1(M_n)$ grows linearly with $n$.
\item[B.] There exists a primitive class $x\in H^1(M;\BZ)$ such that for
$\textbf{any}$ basis $\{x,y\}$ of $H^1(M;\BZ)$,
$\b^n(x,y)=0$ for all $n\ge1$.
\item[C.] There exists a primitive class $x\in H^1(M;\BZ)$ such that for
$\textbf{some}$ basis $\{x,y\}$ of $H^1(M;\BZ)$, $\b^n(x,y)$ can be defined and contains $0$
for each $n\ge1$.
\end{enumerate}
\end{thm}


\begin{cor}\label{systole} Let $M$ be a closed oriented $3$-manifold with
$\b_1(M)=2$. Let $\wt M$ denote the universal torsion-free abelian ($\textbf{Z}\oplus \textbf{Z}$) cover of $M$. Let $[F]$ denote the class in $H_1(\wt M)$ of a lift of a typical fiber of the Abel-Jacobi map of $M$ (represented by a lift of the circle we called $c(x,y)$ below). If, for \textbf{some} $\{x,y\}$, and \textbf{some} $n$, $\b^n(x,y)\neq 0$ then $[F]$ is non-zero.
\end{cor}

The above Corollary generalizes an (independent) result of A. Marin (see Prop.12.1 of [KL]), which dealt with only the case $n=1$. The significance of this Corollary is that it has been previously shown by Ivanov and Katz  ([IK, Theorem 9.2 and Cor.9.3]) that the conclusion of Corollary~\ref{systole} is sufficient to guarantee a certain optimal systolic inequality for $M$. The interested reader is referred to those works.

\begin{figure}[htbp]
\setlength{\unitlength}{1pt}
\begin{picture}(105,92)
\put(10,10){\includegraphics{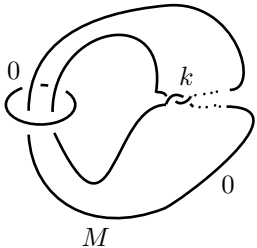}}
\put(11,63){$0$}
\put(92,20){$0$}
\put(76,61){$k$}
\put(39,0){$M$}
\end{picture}
\caption{Example of $\b^1(x,y)= -k$} \label{satolevine}
\end{figure}

Suppose $c$ and $d$ are disjointly embedded oriented circles in $M$ that
are zero in $H_1(M;\BQ)$. Then the {\bf linking number of $c$ with $d$},
$\lk(c,d)\in\BQ$ is defined as follows. Choose an embedded oriented
surface $V_d$ whose boundary is ``$m$ times $d$'' (i.e. a circle in a
regular neighborhood $N$ of $d$ that is homotopic in $N$ to $md$) for some
positive integer $m$, and set:
$$
\lk(c,d) = \f1m(V_d\cd c).
$$
Given this, the invariants $\b^n(x,y)$ are defined as
follows. Let $\{V_x,V_y\}$ be embedded, oriented connected surfaces that
are Poincar\'e Dual to $\{x,y\}$ and meet transversely in an oriented
circle that we call $c(x,y)$ (by the proof of [C1, Theorem 4.1] we may
assume that $c(x,y)$ is connected). Let $c^+(x,y)$ denote a parallel of
$c(x,y)$ in the direction given by $V_y$. Note that $\{V_x,V_y\}$ induce
two maps $\psi_x$, $\psi_y$ from $M$ to $S^1$ wherein the surfaces arise
as inverse images of a regular value. The product of these maps yields a
map $\psi:M\to S^1\x S^1$ that induces an isomorphism on $H_1$/torsion.
Since $c(x,y)$ and $c^+(x,y)$ are mapped to points under $\psi$, they
represent the zero class in $H_1(M;\BQ)$. Therefore we may define
$\b^1(x,y)=\lk(c(x,y)$, $c^+(x,y))$. In fact, $-\b^1(x,y)\cd|\tor
H_1(M;\BZ)|$ is precisely Lescop's invariant of $M$ [Les; p.90-94]. An example is shown in Figure~\ref{satolevine} of a manifold with $\b^1(x,y)= -k$.

The idea of the higher invariants is to iterate this process as long as
possible (compare [C1]). Since $c^+(x,y)$ is rationally null-homologous, there is a surface $V_{c(x,y)}$ whose boundary is ``$k$ times
$c^+(x,y)$'' (in the sense above). We could then define $c(x,x,y)$ to be
$V_x\cap V_{c(x,y)}$, an embedded oriented circle on $V_x$. If $c(x,x,y)$
is rationally null-homologous, then $\b^2(x,y)$ is defined as
$\lk(c(x,x,y),c^+(x,x,y))$ and we may also continue and define
$c(x,x,x,y)$. In general $c(x,x\dots,x,y)=c(x^n,y)$ will be able to be
defined using the chosen surfaces if $c(x^{n-1},y)$ is defined and is also rationally null-homologous (but to do so involves one more choice of a bounding surface). Once $c(x^{n},y)$ is defined and is rationally null-homologous, we may define $\b^n(x,y)$. In general, we do not claim that the value of $\b^n(x,y)$ is independent of the choices of surfaces. Therefore the invariants can be thought of as taking values in a set, just like Massey products. This indeterminacy will not concern us here, for we are only interested in the first non-vanishing value (if it exists) and we shall see that this is independent of the surfaces. 

Much of the time it is convenient to abbreviate $c(\overbrace{x\dots x}^n,y)$ as $c(n)$ so $c(x,y)=c(1)$.

\begin{defn} If $c(n)$ is defined and rationally null-homologous then
$\b^n(x,y)$ is defined to be the set of rational numbers $\lk(c(n),c^+(n))$, ranging over all possible ways of defining such a $c(n)$. If no such $c(n)$ exists then $\b^n(x,y)$ is undefined.
\end{defn}

\begin{example}\label{Example3}
Consider the manifold $M$, shown in Figure~\ref{example3}, obtained from $0$-framed surgery on a two component link $\{ L_x,L_y\}$. Use a genus one Seifert surface for $L_y$ obtained from the obvious twice-punctured disk and a tube that goes up to avoid $L_x$. Let $V_y$ be this surface capped-off in $M$. Similarly use the fairly obvious Seifert surface for $L_x$ in the complement of $L_y$. Then $c^+(x,y)$ is shown. Since it has self-linking zero with respect to $V_x$, $\b^1(x,y)=\b^1(y,x)=0$. Furthermore $\b^2(x,y)=-1$ (note the link $\{c^+(x,y),L_x\}$ is very similar to that of Figure~\ref{satolevine}). This means that $\pi_1(M)$ does \textbf{not} admit an epimorphism to $\mathbb{Z}\ast \mathbb{Z}$ since that would imply that $\{ L_x,L_y\}$ were a homology boundary link. But $\b^2(x,y)=-1$ precludes this by \cite{C2}. Nonetheless, further $c(yy...y,x)$ may be taken to be empty since $c^+(x,y)$ and $L_y$ form a boundary link in the complement of $L_x$. Thus $\b^n(y,x)= 0$ for all $n$, indicating, by Theorem~\ref{linear}, that the first Betti numbers will grow linearly in the family of finite cyclic covers corresponding to the map $\pi_1(M)\to \mathbb{Z}$ that sends a meridian of $L_x$ to zero and a meridian of $L_y$ to one.
\begin{figure}[htbp]
\setlength{\unitlength}{1pt}
\begin{picture}(123,73)
\put(10,10){\includegraphics{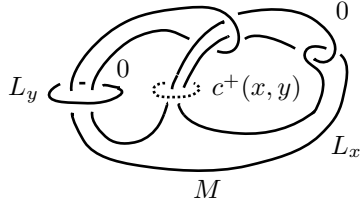}}
\put(72,39){$c^+(x,y)$}
\put(-5,39){$L_y$}
\put(117,17){$L_x$}
\put(65,0){$M$}
\put(36,46){$0$}
\put(119,68){$0$}
\end{picture}
\caption{Example with linear growth in cyclic covers but no map to $\mathbb{Z}\ast\mathbb{Z}$} \label{example3}
\end{figure}
\end{example}
\begin{example}\label{Example2}
Consider the family of manifolds $M_k$, shown in Figure~\ref{example2} and Figure~\ref{example2b}, obtained from $0$-framed surgery on a two component link. 
\begin{figure}[htbp]
\setlength{\unitlength}{1pt}
\begin{picture}(135,98)
\put(10,10){\includegraphics{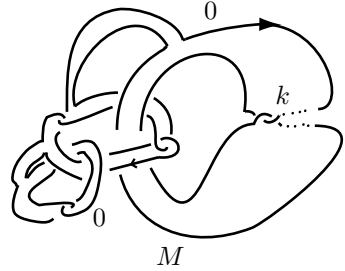}}
\put(83,93){$0$}
\put(41,15){$0$}
\put(110,61){$k$}
\put(65,0){$M$}
\end{picture}
\caption{Example with $\b^1(x,y)=0$,$\b^2(x,y)=-k$, $\b^2(y,x)= -1$} \label{example2}
\end{figure}
\begin{figure}[htbp]
\setlength{\unitlength}{1pt}
\begin{picture}(135,98)
\put(10,10){\includegraphics{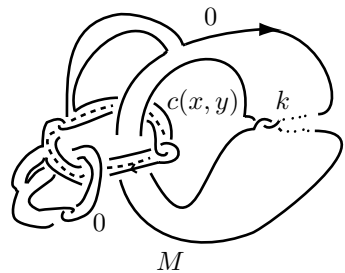}}
\put(83,93){$0$}
\put(41,15){$0$}
\put(110,61){$k$}
\put(69,61){$c(x,y)$}
\put(65,0){$M$}
\end{picture}
\caption{The circle $c(x,y)$} \label{example2b}
\end{figure}

If $V_x$ denotes the capped-off Seifert surface (obtained using Seifert's algorithm) for the link component, $L_x$, on the right-hand side and $V_y$ denotes the capped-off Seifert surface for the link component, $L_y$, on the left-hand side, then the dashed circle in Figure~\ref{example2b} is $c(x,y)=V_x\cap V_y$. The circle $c(y,x)$ is merely this circle with opposite orientation. Since it lies on an untwisted band of $V_x$, $\b^1(x,y)=0=\b^1(y,x)$. Therefore the Lescop invariant of $M$ vanishes. But the link $\{c(x,y),L_x\}$ is the link of Figure~\ref{satolevine} so $\b^2(x,y)=-k$, whereas the link $\{c(y,x),L_y\}$ is a Whitehead link so $\b^2(y,x)= -1$. We claim further that, as long as $k\neq 0$,  for $\textbf{any}$ basis $\{X,Y\}$ of $H^1(M)$, $\b^2(X,Y)\neq 0$. It will then follow from Theorem~\ref{linear} that the first Betti number of $M$ will grow sub-linearly in \textbf{any} family of finite cyclic covers. A general basis, $\{V_X,V_Y\}$, of $H_2(M)$ can be represented as follows. Represent $V_X$ by $p$ parallel copies of $V_x$ together with $q$ parallel copies $V_y$, and represent $V_Y$ by $r$ parallel copies of $V_x$ together with $s$ parallel copies $V_y$, where $ps-qr=\pm 1$. Thus $c(X,Y)=-c(Y,X)$ is represented by $ps-qr$ parallel copies of $c(x,y)$. It follows that $\b^1(X,Y)=\b^1(x,y)=0$, reinforcing our above claim that $\b^1$ is independent of basis. Hence $V_{c(X,Y)}=\pm V_{c(x,y)}$ so $c(X,X,Y)$ is represented by $\pm pc(x,x,y)\mp qc(y,y,x)$. Since $\b^2(X,Y)$ is the self-linking number of this class, it can be evaluated to be
$$
p^2\b^2(x,y) + q^2\b^2(y,x) - 2pq\lk (c(x,x,y),c(y,y,x))
$$
but the latter mixed linking number is easily seen to be zero in this case. Hence $\b^2(X,Y)=-kp^2-q^2$ which is non-zero if $k$ is non-zero.
\end{example}

\begin{lem}\label{equivalence} Suppose $c(1),\dots,c(n)$ have been defined
as embedded oriented curves on $V_x$ arising as $c(1)=V_x\cap V_y$ and
$$
c(j) = V_x\cap V_{c(j-1)}\qquad2\le j\le n
$$
where $V_{c(j)}$, $1\le j\le n-1$, is an embedded, oriented connected
surface whose boundary is a positive multiple $k_j$ of $c^+(j)$ (in the
sense above). Then $\b^j$ is defined for $1\le j\le n-1$ and the following are equivalent:
\begin{enumerate} 
\item[B1:] $\b^1,\dots,\b^n$ are defined using the given system of surfaces.
\item[B2:] $\b^j$ is defined for $1\le j\le n$ and is {\bf zero} for $1\le
j\le\[\f n2\]$
\item[B3:] $c(n+1)$ exists
\item[B4:] For all $s$, $t$ such that $1\le s\le t$ and $s+t\le n$
, $\lk(c(s),c^+(t))=0$.
\end{enumerate}
\end{lem}

\begin{proof}[Proof of Proposition~\ref{equivalence}] Assume $1\le j\le
n-1$. The hypotheses imply that a positive multiple of $c^+(j)$ is
(homotopic to) the boundary of a surface so $c^+(j)$ and $c(j)$ are
rationally null-homologous. Thus their linking number is well-defined,
establishing the first claim.

\medbreak
\noindent{\bf B1$\Longleftrightarrow$B3}: $\b^n$ is defined precisely when
$[c(n)]=0$ in $H_1(M;\BQ)$ which is precisely the condition under which
$c(n+1)$ can be defined.

\medbreak
\noindent{\bf B1$\Longrightarrow$B4}: If $n=1$ the implication is true
since B4 is vacuous. Thus assume by induction that the implication is true
for $n-1$, that is our inductive assumption is that, for all $s+t<n$, $\lk(c(s),c^+(t))=0$. Now  consider the case
that $s+t=n$. Since $\b^n$ is defined $[c(n)]=0$ in $H_1(M;\BQ)$. We
claim this is true precisely when $c(n)\cd c(1)=0$ (here we refer to oriented intersection number on the surface $V_x$). For suppose
$\psi_x:M\to S^1$ and $\psi_y:M\to S^1$ are maps such that
$\psi^{-1}_x(*)=V_x$ and $\psi^{-1}_y(*)=V_y$. Then
$(\psi_x)_*([c(n)])=0$ since $c(n)\subset V_x$; and $(\psi_y)_*([c(n)])=0$
precisely when $c(n)\cd V_y=c(n)\cd c(1)=0$. But the map $\psi_x \times \psi_y$ completely
detects $H_1(M)$/Torsion. Therefore, once $c(n)$ exists, $\b^n$ is defined if and only if:
\begin{align*}
0 = c(n)\cd c(1)  &= \pm(c(1)\cd V_{c(n-1)})\\
&= \pm k_{n-1}\lk(c(1),c^+(n-1))
\end{align*}
which establishes B4 in the case $s=1$. But we claim that, if B4 is true
for $s+t<n$, then for $s+t=n$ and $s< t$,
$$
k_{t-1}\lk(c(s+1),c^+(t-1))=k_s\lk(c(s),c^+(t)).
$$
This equality can then be applied, successively decreasing $s$, to establish B4 in generality. This claimed equality
is established as follows. 
\begin{align*}
\pm k_{t-1}\lk(c(s+1),c^+(t-1)) &= \pm V_{c(t-1)}\cd c(s+1)\\
&= \pm V_{c(t-1)}\cd(V_x\cap V_{c(s)})\\
&= V_{c(s)}\cd(V_x\cap V_{c(t-1)})\\
&= V_{c(s)}\cd c(t)\\
&= k_s\lk(c(t),c^+(s))\\
&=k_s\lk(c(s),c^+(t)).
\end{align*}
The last step is justified by verifying that $c(s)\cd c(t)=0$ if $s< t$.
For
\begin{align*}
c(s)\cd c(t)  &= c(s)\cd V_{c(t-1)}\\
&= \pm k_{t-1}\lk(c(s),c^+(t-1))
\end{align*}
which vanishes by our inductive assumption since $s+(t-1)<n$. 

\medbreak
\noindent{\bf B4$\Longrightarrow$B1}: Since $\b^1$ is always defined we
may assume $n>1$. It follows from B4 that $\lk(c(1),c^+(n-1))=0$ if $n>1$.
But we saw in the proof of B1$\Longrightarrow$B4 that once $c(n)$ was defined, this was equivalent
to $\b^n$ being defined.

\medbreak
\noindent{\bf B2$\Longrightarrow$B1}: This is obvious.

\medbreak
\noindent{\bf B4$\Longrightarrow$B2}: Since B4$\Longrightarrow$B1, we have
$\b^j$ defined for $j\le n$. Now suppose $1\le j\le\[\f n2\]$. Since
$\b^j=\lk(c(j),c^+(j)$ and $2j\le n$, this vanishes by B4.

This completes the proof of Lemma~\ref{equivalence}.
\end{proof}

\begin{proof}[Proof of Theorem~\ref{linear}] The proof shows slightly more, namely that there is a correspondence between the infinite cyclic cover implicit in part $A$ and the class $x$ in parts $B$ and $C$. Suppose $\{M_n\}$ is a family
of $n$-fold cyclic covers of $M$ corresponding to the infinite cyclic
cover $M_\infty$. Note that $H_1(M_\infty;\BQ)$ is a finitely generated
$\La=\BQ[t,t^{-1}]$ module (this involves a choice of generator of the
infinite cyclic group of deck translations of $M_\infty$). Throughout this proof, homology will be taken with rational coefficients unless specified otherwise.

\medbreak
\noindent{\bf Step 1}: $\b_1(M_n)$ grows linearly
$\Longleftrightarrow H_1(M_\infty;\BQ)$ has positive rank as a $\La$-module.

\medbreak
As remarked above, this fact was previously known. We present a quick proof for
the convenience of the reader. We are indebted to Shelly Harvey for
showing us this elementary proof. Since $\La$ is a PID,
$$
H_1(M_\infty)\cong\La^{r_1}\oplus_j\f\La{\<p_j(t)\>}
$$
where $p_j(t)\neq0$. By examining the ``Wang sequence'' with $\BQ$-coefficients
$$
H_2(M_\infty)\lra H_2(M_n)\overset{\p_*}{\lra}H_1(M_\infty)
\overset{t^n-1}{\lra}H_1(M_\infty)\overset{\pi}{\lra}H_1(M_n)
\overset{\p_*}{\lra}H_0(M_\infty)\overset{t^n-1}{\lra}
$$
 it is easily seen that
\begin{align*}
H_1(M_n)  &\cong\f{H_1(M_\infty)}{\<t^n-1\>}\op\BQ\\
&\cong\(\f\La{\<t^n-1\>}\)^{r_1}\oplus_j\f\La{\<p_j(t),t^n-1\>}\op\BQ.
\end{align*}
The first summand contributes $nr_1$ to $\b_1(M_n)$. The $\BQ$-rank of the
second summand is bounded above by the sum of the degrees of the $p_j$, a
number that is {\bf independent} of $n$. Therefore $\b_1(M_n)$ grows
linearly with $n$ if $r_1\neq0$ and otherwise is bounded above by a
constant (independent of $n$).

\medbreak
\noindent{\bf Step 2}: $H_1(M_\infty)$ has positive $\La$-rank
$\Longleftrightarrow H_1(M_\infty)$ has no $(t-1)$-torsion (equivalently $t-1$ acts injectively).
To verify Step 2, consider the ``Wang sequence'' with $\BQ$-coefficients
$$
H_2(M_\infty)\lra H_2(M)\overset{\p_*}{\lra}H_1(M_\infty)
\overset{t-1}{\lra}H_1(M_\infty)\overset{\pi}{\lra}H_1(M)
\overset{\p_*}{\lra}H_0(M_\infty)\overset{t-1}{\lra}H_0(M_\infty)
$$
associated to the exact sequence of chain complexes
$$
0\lra C_*(M_\infty;\BQ)\overset{t-1}{\lra}
C_*(M_\infty;\BQ)\overset{\pi}{\lra}C_*(M;\BQ)\lra0.
$$
Since $H_0(M_\infty)\cong\BQ$, $\image\p_*\cong\BQ$ on $H_1(M)$. If
$\b_1(M)=2$ then it follows that
$\BQ\cong\ke\p_*=\image(\pi)\cong\cok(t-1)$. It follows that $H_1(M_\infty)$
contains at most one summand of the form $\La/\<(t-1)^m\>$ since each such
summand contributes precisely one $\BQ$ to $\cok(t-1)$. Similarly each
$\La$ summand of $H_1(M_\infty)$ contributes one $\BQ$ to the cokernel.
Therefore $H_1(M_\infty)$ has positive $\La$ rank if and only if it has no
summand of the form $\La/\<(t-1)^m\>$. The latter is equivalent to saying
that it has no $(t-1)$-torsion, or that $t-1$ acting on $H_1(M_\infty)$ is
injective. This completes Step 2.

\medbreak
\noindent{\bf Step 3}: $(t-1):H_1(M_\infty)\to H_1(M_\infty)$ is injective
$\Longleftrightarrow$ For any surface $V_x$, dual to $x$, and for each surface $V_y$ such that $\{[V_y],[V_x]\}$ generates $H_2(M;\BZ)$, the class $[\tl c(x,y)]\in
H_1(M_\infty;\BQ)$ is zero. Moreover the latter statement is equivalent to
one where ``for each'' is replaced by ``for some''.

Suppose that $t-1$ is
injective. Note that the injectivity of $t-1$ is equivalent to
$\p_*:H_2(M)\to H_1(M_\infty)$ being the zero map. Then, for {\bf any} $[V_y]$ as above, $\p_*([V_y])=0$. But we
claim that $\p_*([V_y])$ is represented by $[\tl c(x,y)]$, since $V_x$ is
Poincar\'e Dual to the class $x$ defining $M_\infty$. For if $Y=M-\intt(V_x\x[-1,1])$  then a copy of $Y$, denoted $\wt Y$, can be
viewed as a fundamental domain in $M_\infty$, as shown in Figure~\ref{cover}.
\begin{figure}[htbp]
\setlength{\unitlength}{1pt}
\begin{picture}(177,84)
\put(10,10){\includegraphics{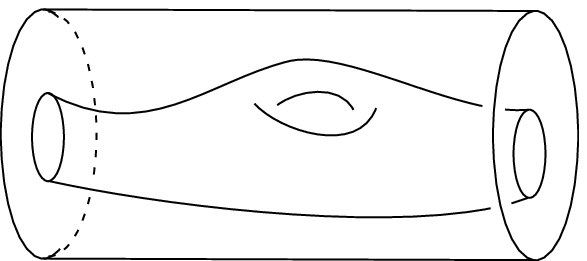}}
\put(15,62){$\wt c(x,y)$}
\put(80,33){$\wt V_y$}
\put(0,12){$\wt V_x$}
\put(173,12){$t\wt V_x$}
\put(179,59){$t\wt c(x,y)$}
\put(87,-5){$\wt Y$}
\end{picture}
\caption{Fundamental Domain of $M_\infty$} \label{cover}
\end{figure}

Moreover if $\wt V_y$ denotes $p^{-1}(V_y)\cap \wt Y$ then $\wt V_y$ is a compact surface in $\wt M$ whose boundary is $t_*(\tl c(x,y))-\tl c(x,y)$. Thus $\wt V_y$ is a 2-chain in
$M_\infty$ such that $\pi_\#(\wt V_y)$ gives the chain representing
$[V_y]$. Since $\p\wt V_y$ is $(t-1)\tl c(x,y)$ in $C_*(M_\infty;\BQ)$, it follows from the explicit construction of
$\p_*$ in the proof of the Zig-Zag Lemma [Mu,Section 24] that $\p_*([V_y])=[\tl
c(x,y)]$. 

Conversely, if $\p_*([V_y])=0$ for {\bf some} $[V_y]$ then
$\p_*$ is the zero map (note that since $V_x$ lifts to $M_\infty$, $[V_x]$ lies in the
image of $H_2(M_\infty)\lra H_2(M)$ so $\p_*([V_x])=0$ t).

Therefore the injectivity statement implies the ``for each'' statement which
clearly implies the ``for some'' statement. Conversely, the ``for some''
statement implies the injectivity statement.

\medbreak
\noindent{\bf Step 4}: The class $[\tl c(x,y)]$ from Step 3 is 0 if and
only if it is divisible by $(t-1)^k$ for every positive $k$. In fact it suffices that it be divisible by $(t-1)^N$ where $N$ is the largest nonnegative integer such that $\La/\<(t-1)^N\>$ is a summand of $H_1(M_\infty,\BQ)$.

\medbreak
One implication is immediate, so assume that there exists a class $[V_y]$
as in Step 3 such that $\p_*([V_y])=[\tl c_{1}]=(t-1)^N\b$ for some $\b\in H_1(M_\infty)$. Since $[\tl c_{1}]\in\image\p_*$, it is $(t-1)$-torsion so $\b$ is $(t-1)^{N+1}$-torsion. Moreover $\b$ lies in the submodule $A\subset
H_1(M_\infty,\BQ)$ consisting of elements annihilated by some power of
$t-1$, so, by choice of $N$, $(t-1)^N\b=0=[\tl c_{1}]=0$ as desired.
This completes the verification of Step 4.

\medbreak
\noindent{\bf Step 5}: C$\Rightarrow$A 
Let $\{x,y\}$ be as in the
hypotheses of C and let $M_\infty$ correspond to the class $x$. Let $N$ be the positive integer as above. If $\b^{(N+1}$ can be defined, we know in particular that there exists some system of surfaces $\{V_x,V_y,...,V_{c(N)}\}$ that defines $\{c(j)\}$, $1\leq j\leq (N+1)$. Choose
a preferred lift $\wt V_x$, of $V_x$ to $M_\infty$ and a preferred fundamental domain $\wt Y$ as above lying on the positive side of $\wt V_x$. Consider any $m, 1\leq m \leq N$. Since $c(m)$ and $c(m+1)$ lie on $V_x$, they lift to oriented 1-manifolds
$\tl c(m)$ and $\tl c(m+1)$ in $\wt V_x$. Similarly $c^+(m)$ lifts to $\tl c^+(m)$, which is a push-off of $\tl c(m)$ lying in $\wt Y$. Recall that
$c(m+1)=V_{c(m)}\cap V_x$ where $\p
V_{c(m)}=k_mc^+_{(m)}$ for some positive integer $k_m$. Letting
$\wt V_{c(m)}$ be $V_{c(m)}$ cut open along $c(m+1)$ we
observe that $\wt V_{c(m)}$ can be lifted to $\wt Y$ and viewed as a 2-chain
showing that $k_m[\tl c^+(m)]=(t-1)[\tl c(m+1)]$ in
$H_1(M_\infty;\BQ)$, as in Figure~\ref{stepfive}. 

\begin{figure}[htbp]
\setlength{\unitlength}{1pt}
\begin{picture}(177,84)
\put(10,10){\includegraphics{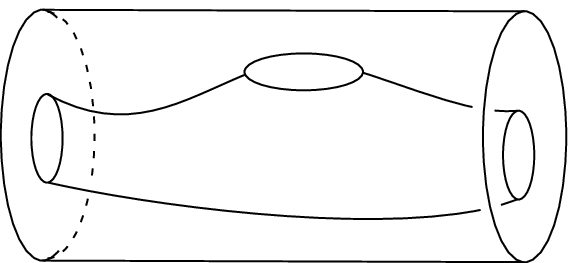}}
\put(15,64){$\wt c(m+1)$}
\put(108,72){$k_m\wt c^+(m)$}
\put(65,34){$\wt V_{c(m)}$}
\put(0,12){$\wt V_x$}
\put(173,12){$t\wt V_x$}
\put(179,59){$t\wt c(m+1)$}
\put(83,-5){$\wt Y$}

\end{picture}
\caption{} \label{stepfive}
\end{figure}

Thus $[\tl c^+(1)]=(t-1)(1/k_1)[\tl c(2)]=(t-1)^2(1/k_1)(1/k_2)[\tl c(3)]$, et cetera, showing that $[\tl c(1)]$ is divisible by $(t-1)^N$. By Steps 1 through 4, this implies A of
Theorem~\ref{linear}, completing Step 5.
\medbreak
\noindent{\bf Step 6}: A$\Rightarrow$B We assume that there is a primitive
class $x\in H^1(M;\BZ)$ corresponding to $M_\infty$ and $\{M_n\}$ where
$\b_1(M_n)$ grows linearly. By Steps 1, 2, and 3, for any $\{x,y\}$
generating $H^1(M;\BZ)$,  $H_1(M_\infty;\BQ)$ has no
$(t-1)$-torsion, and for any surfaces dual to ${x,y}$, $[\tl c(1)]=0$. Recall that $c(1)$ and $c(2)$ are always defined.
We shall establish inductively that for all $m\ge2$, $c(m)$ is defined
and that for \emph{any} system of surfaces used to define $c(m-1)$, $[\tl c(m-1)]=0$ in $H_1(M_\infty;\BQ)$. This has already been
shown for $m=2$. Suppose it has been established for $m$ (and all lesser
values). We now establish it for $m+1$. Since $c(m)$ and
$c(m-1)$ exist, the argument in Step 5 (Figure~\ref{stepfive}) shows that
 $k_{m-1}[\tl c(m-1)]=(t-1)[\tl c(m)]$ in $H_1(M_\infty;\BQ)$. But $[\tl c(m-1)]=0$ so $[\tl c(m)]$ is $(t-1)$-torsion. Since there is no
non-trivial $(t-1)$-torsion, $[\tl c(m)]=0$ in $H_1(M_\infty;\BQ)$.
Hence $[c(m)]=0$ in $H_1(M;\BQ)$ so $c(m+1)$ is defined.
Since this holds for any system of defining surfaces, this completes the
inductive step.

Since $c(m)$ is defined for all $m\ge1$, by Lemma~\ref{equivalence},
$\b^n(x,y)=0$ for all $n\ge1$. This completes the proof of Step 6.

Since B clearly implies C, this completes the proof of
Theorem~\ref{linear}.
\end{proof}
\begin{proof}[Proof of Corollary~\ref{systole}] Assume some $\b^m(x,y)\neq 0$. If $[F]$ were zero then certainly, for the fixed infinite cyclic cover, $M_\infty$, corresponding to $x$, $[\wt c(x,y)]=0$ in $H_1(M_\infty;\mathbb{Q})$ so by Steps $1-3$ of the above proof, $\b_1(M_n)$ grows linearly with $n$. By Theorem~\ref{linear}, this would imply that $\b^m(x,y)= 0$ for all $m$, a contradiction.
\end{proof}

\end{document}

%% file: amsmath-defs.tex
\def\b{\beta}
\def\cd{\cdot}

\def\f{\frac}

\def\lb\{{\left\{}
\def\la{\lambda}
\def\La{\Lambda}
\def\lla{\longleftarrow}
\def\lm{\limits}
\def\lra{\longrightarrow}
\def\dllra{\Longleftrightarrow}
\def\llra{\longleftrightarrow}
\def\n{\nabla}
\def\ngth{\negthickspace}
\def\ola{\overleftarrow}
\def\Om{\Omega}
\def\om{\omega}
\def\op{\oplus}
\def\oper{\operatorname}
\def\oplm{\operatornamewithlimits}
\def\ora{\overrightarrow}
\def\ov{\overline}
\def\ova{\overarrow}
\def\ox{\otimes}
\def\p{\partial}
\def\rb\}{\right\}}
\def\s{\sigma}
\def\sbq{\subseteq}
\def\spq{\supseteq}
\def\sqp{\sqsupset}
\def\supth{{\text{th}}}
\def\T{\Theta}
\def\th{\theta}
\def\tl{\tilde}
\def\thra{\twoheadrightarrow}
\def\un{\underline}
\def\ups{\upsilon}
\def\vp{\varphi}
\def\wh{\widehat}
\def\wt{\widetilde}
\def\x{\times}
\def\z{\zeta}
\def\({\left(}
\def\){\right)}
\def\[{\left[}
\def\]{\right]}
\def\<{\left<}
\def\>{\right>}

\def\tec{Teichm\"uller\ }
\def\sconr{\hbox{\medspace\vrule width 0.4pt height 4.7pt depth
0.4pt \vrule width 5pt height 0pt depth 0.4pt\medspace}}

\def\SA{\mathcal A}
\def\SB{\mathcal B}
\def\SC{\mathcal C}
\def\SD{\mathcal D}
\def\SE{\mathcal E}
\def\SF{\mathcal F}
\def\SG{\mathcal G}
\def\SH{\mathcal H}
\def\SI{\mathcal I}
\def\SJ{\mathcal J}
\def\SK{\mathcal K}
\def\SL{\mathcal L}
\def\SM{\mathcal M}
\def\SN{\mathcal N}
\def\SO{\mathcal O}
\def\SP{\mathcal P}
\def\SQ{\mathcal Q}
\def\SR{\mathcal R}
\def\SS{\mathcal S}
\def\ST{\mathcal T}
\def\SU{\mathcal U}
\def\SV{\mathcal V}
\def\SW{\mathcal W}
\def\SX{\mathcal X}
\def\SY{\mathcal Y}
\def\SZ{\mathcal Z}

\newcommand{\BA}{\ensuremath{\mathbf A}}
\newcommand{\BB}{\ensuremath{\mathbf B}}
\newcommand{\BC}{\ensuremath{\mathbf C}}
\newcommand{\BD}{\ensuremath{\mathbf D}}
\newcommand{\BE}{\ensuremath{\mathbf E}}
\newcommand{\BF}{\ensuremath{\mathbf F}}
\newcommand{\BG}{\ensuremath{\mathbf G}}
\newcommand{\BH}{\ensuremath{\mathbf H}}
\newcommand{\BI}{\ensuremath{\mathbf I}}
\newcommand{\BJ}{\ensuremath{\mathbf J}}
\newcommand{\BK}{\ensuremath{\mathbf K}}
\newcommand{\BL}{\ensuremath{\mathbf L}}
\newcommand{\BM}{\ensuremath{\mathbf M}}
\newcommand{\BN}{\ensuremath{\mathbf N}}
\newcommand{\BO}{\ensuremath{\mathbf O}}
\newcommand{\BP}{\ensuremath{\mathbf P}}
\newcommand{\BQ}{\ensuremath{\mathbf Q}}
\newcommand{\BR}{\ensuremath{\mathbf R}}
\newcommand{\BS}{\ensuremath{\mathbf S}}
\newcommand{\BT}{\ensuremath{\mathbf T}}
\newcommand{\BU}{\ensuremath{\mathbf U}}
\newcommand{\BV}{\ensuremath{\mathbf V}}
\newcommand{\BW}{\ensuremath{\mathbf W}}
\newcommand{\BX}{\ensuremath{\mathbf X}}
\newcommand{\BY}{\ensuremath{\mathbf Y}}
\newcommand{\BZ}{\ensuremath{\mathbf Z}}

\def\bba{{\mathbb A}}
\def\bbb{{\mathbb B}}
\def\bbc{{\mathbb C}}
\def\bbd{{\mathbb D}}
\def\bbe{{\mathbb E}}
\def\bbf{{\mathbb F}}
\def\bbg{{\mathbb G}}
\def\bbh{{\mathbb H}}
\def\bbi{{\mathbb I}}
\def\bbj{{\mathbb J}}
\def\bbk{{\mathbb K}}
\def\bbl{{\mathbb L}}
\def\bbm{{\mathbb M}}
\def\bbn{{\mathbb N}}
\def\bbo{{\mathbb O}}
\def\bbp{{\mathbb P}}
\def\bbq{{\mathbb Q}}
\def\bbr{{\mathbb R}}
\def\bbs{{\mathbb S}}
\def\bbt{{\mathbb T}}
\def\bbu{{\mathbb U}}
\def\bbv{{\mathbb V}}
\def\bbw{{\mathbb W}}
\def\bbx{{\mathbb X}}
\def\bby{{\mathbb Y}}
\def\bbz{{\mathbb Z}}